\documentclass[a4paper,12pt]{article}
\begin{document}

\title{Cohomologies of nonassociative metagroup algebras.}
\author{S.V. Ludkowski.}

\date{25 June 2017}
\maketitle

\begin{abstract}
Cohomologies of nonassociative metagroup algebras are investigated.
Extensions of metagroup algebras are studied. Examples are given.
\footnote{key words: nonassociative algebra; metagroup; cohomology; extension 
\\ address: Dep. Appl. Mathematics, Moscow State Techn. Univ. MIREA,
av. Vernadsky 78, Moscow, 119454, Russia, \par e-mail:
sludkowski@mail.ru }

\end{abstract}

\section{Introduction.}
\par Nonassociative algebras compose a very important area of algebra.
Among them nonassociative Cayley-Dickson algebras are intensively
studied. They are generalizations of the octonion algebra (see
\cite{baez,dickson,kansol,schaeferb} and references therein). Their
structure and identities in them attract great attention (see, for
example, \cite{baez,girardgb,guertzeb,kansol,krausshb,serodaaca07}
and references therein). They have found many-sided applications in
the theory of Lie groups and algebras and their generalizations,
mathematical analysis, non-commutative geometry, operator theory,
PDE and their applications in natural sciences including physics and
quantum field theory. \par An extensive area of investigations of
PDE intersects with cohomologies and deformed cohomologies
\cite{pommb}. Therefore, it is important to develop this area over
octonions, Cayley-Dickson algebras and more general metagroup
algebras.
\par It appears that generators of
Cayley-Dickson algebras form objects, which are nonassociative
generalizations of groups. They are called metagroups. This means
that metagroup algebras include as the particular case the
Cayley-Dickson algebras. This article is devoted to algebras
generated by metagroups.
\par On the other hand, algebras are frequently studied using
cohomology theory. But the already developed cohomology theory
operates with associative algebras. It was investigated by
Hochschild and other authors
\cite{bourbalgbhomol,bredonb67,cartaneilenbb56,hochschild46}, but it
is not applicable to nonassociative algebras. This work is devoted
to the development of cohomology theory for nonassociative algebras,
namely for its subclass of metagroup algebras.
\par Previously cohomologies of loop spaces on quaternion and octonion manifolds
were studied in \cite{ludwrgrijmgta}. They have specific features in
comparison with the case of complex manifolds. This is especially
caused by the noncommutativity of the quaternion skew field and the
nonassociativity of the octonion algebra.
\par In this article metagroup algebras which may be nonassociative
are studied. Their cohomologies are investigated. For this purpose
two-sided modules also are considered. Extensions of metagroup
algebras are investigated. The corresponding cohomology theory is
developed for this purpose. Moreover, homomorphisms and
automorphisms of metagroup algebras are investigated.
\par All main results of this paper are obtained for the first time.
They can be used for further studies of nonassociative algebra
cohomologies, structure of nonassociative algebras, operator theory
and spectral theory over Cayley-Dickson algebras, PDE,
non-commutative analysis, non-commutative geometry, mathematical
physics, their applications in the sciences.

\section{Cohomologies of nonassociative metagroup algebras.}

\par To avoid misunderstandings we give necessary definitions.

\par {\bf 1. Definitions.}  Let $G$ be a set with a single-valued
binary operation (multiplication)  $G^2\ni (a,b)\mapsto ab \in G$
defined on $G$ satisfying the conditions: \par $(1)$ for each $a$
and $b$ in $G$ there is a unique $x\in G$ with $ax=b$ and \par $(2)$
a unique $y\in G$ exists satisfying $ya=b$, which are denoted by the
left division $x=a\setminus b$ and the right division $y=b/a$
correspondingly,
\par $(3)$ there exists a neutral (i.e. unit) element $e\in G$: $~eg=ge=g$
for each $g\in G$. In another notation it also is used $1$ instead
of $e$.
\par The set of all elements $h\in G$
commuting and associating with $G$:
\par $(4)$ $Com (G) := \{ a\in G: \forall b\in G ~ ab=ba \} $,
\par $(5)$ $N_l(G) := \{a\in G: \forall b\in G \forall c\in G ~ (ab)c=a(bc) \}
$,
\par $(6)$ $N_m(G) := \{a\in G: \forall b\in G \forall c\in G ~ (ba)c=b(ac)
\} $,
\par $(7)$ $N_r(G) := \{a\in G: \forall b\in G \forall c\in G ~ (bc)a=b(ca)
\} $,
\par $(8)$ $N(G) := N_l(G)\cap N_m(G)\cap N_r(G)$, $~{\cal C}(G) := Com (G)\cap N(G)$.
\\ is called the center ${\cal C}(G)$ of $G$.
\par Let $F^{\times }$ denote a (proper or improper) subgroup
of ${\cal C}(G)$.
\par We call $G$ a metagroup if a set $G$ possesses a single-valued binary operation
and satisfies conditions $(1-3)$ and
\par $(9)$ $(ab)c={\sf t}_3(a,b,c)a(bc)$ \\ for each
$a$, $b$ and $c$ in $G$, where ${\sf t}_3(a,b,c)\in F^{\times }$, $ ~ F^{\times }\subset {\cal C}(G)$.
\par Then $G$ will be called a central metagroup if in addition to $(9)$ it
satisfies the condition:
\par $(10)$ $ab={\sf t}_2(a,b)ba$ \\ for each $a$ and $b$ in $G$, where
${\sf t}_2(a,b)\in F^{\times }$.
\par In view of the nonassociativity of $G$ in general
a product of several elements of $G$ is specified as usually by
opening "$($" and closing "$)$" parentheses. For elements
$a_1$,...,$a_n$ in $G$ we shall denote shortly by $ \{ a_1,...,a_n
\}_{q(n)}$ the product, where a vector $q(n)$ indicates an order of
pairwise multiplications of elements in the row $a_1,...,a_n$ in
braces in the following manner. Enumerate positions: before $a_1$ by
$1$, between $a_1$ and $a_2$ by $2$,..., by $n$ between $a_{n-1}$
and $a_n$, by $n+1$ after $a_n$. Then put $q_j(n)=(k,m)$ if there
are $k$ opening "$($" and $m$ closing "$)$" parentheses in the
ordered product at the $j$-th position of the type $)...)(...($,
where $k$ and $m$ are nonnegative integers, $q(n)
=(q_1(n),....,q_{n+1}(n))$ with $q_1(n)=(k,0)$ and
$q_{n+1}(n)=(0,m)$.
\par As traditionally $S_n$ denotes the symmetric group of the set $ \{ 1,
2,..., n \} $.

\par {\bf 2. Lemma.} {\it $(i)$. Let $G$ be a central metagroup.
Then for every $a_1$,...,$a_n$ in $G$, $v\in S_n$ and vectors $q(n)$
and $u(n)$ indicating an order of pairwise multiplications and $n\in
\bf N$ there exists an element $t_n=t_n(a_1,...,a_n;q(n),u(n)|v)\in F^{\times
}$ such that
\par $(1)$ $ \{ a_1,...,a_n \}_{q(n)}= t_n \{ a_{v(1)},...,a_{v(n)} \}_{u(n)}$.
\par $(ii)$. If $G$ is a metagroup, then property $(1)$ is satisfied
for the neutral element $v=id$ in $S_n$.}
\par {\bf Proof.} From conditions 1$(1-8)$ it follows that
${\cal C}(G)$ itself is a commutative group. \par $(i)$. For $n=1$
evidently $t_1=1$. For $n=2$ this follows from condition 1$(10)$.
Consider $n=3$. When $u$ is the identity element of $S_3$ the
statement follows from condition 1$(9)$. For any transposition $u$
of two elements of the set $ \{ 1, 2 ,3 \} $ the statement follows
from 1$(9)$ and 1$(10)$. Elements of $S_3$ can be obtained by
multiplication of pairwise transpositions. Therefore from the
condition $F^{\times }\subset {\cal C}(G)$ it follows that formula
$(1)$ is valid.
\par Let now $n\ge 4$ and suppose that this lemma is proved for any
products consisting of less than $n$ elements. In view of properties
1$(1)$ and 1$(2)$ it is sufficient to verify formula $(1)$ for \par
$ \{ a_1,...,a_n \}_{q(n)} = (...((a_1a_2)a_3)...)a_n =: \{
a_1,...,a_n \}_{l(n)}$, \\ since $F^{\times }\subset {\cal C}(G)$.
In the particular case \par $ \{ a_{v(1)},...,a_{v(n)} \}_{u(n)} =
\{ a_{v(1)},...,a_{v(n-1)} \}_{u(n-1)}a_n$ \\ formula $(1)$ follows
from the induction hypothesis, since
\par $(...((a_1a_2)a_3)...)a_{n-1} = t_{n-1} \{
a_{v(1)},...,a_{v(n-1)} \}_{u(n-1)}$ and hence
\par $((...((a_1a_2)a_3)...)a_{n-1})a_n = t_{n-1} (\{
a_{v(1)},...,a_{v(n-1)} \}_{u(n-1)}a_n)$ \\ and putting
$t_n=t_{n-1}$, where $t_{n-1}=t_{n-1}(a_1,...,a_{n-1}; q(n-1), u(n-1)| w)$
with $w=v|_{ \{ 1,...,n-1 \} }$, $ ~ v(n)=n$.
\par In the general case
$ \{ a_{v(1)},...,a_{v(n)} \}_{u(n)}= \{ b_1,...,b_j,...,b_k
\}_{p(k)} $, where $j$ is such that either $b_j = c_ja_n$ with $c_j=\{
a_{v(j)},...,a_{v(j+m-1)} \}_{r(m)}$ and with $v(j+m)=n$ or $b_j =
a_n c_j$ with $c_j=\{ a_{v(j+1)},...,a_{v(j+m)} \}_{r(m)}$ and with
$v(j)=n$, also $b_1=a_{v(1)}$,...,$b_{j-1}=a_{v(j-1)}$,
$b_{j+1}=a_{v(j+1)}$,..., $b_k=a_{v(n)}$ with suitable vectors $p(k)$ and $r(m)$.
If $m>1$, then $k<n$ and
using the induction hypothesis for $\{ b_1,...,b_j,...,b_k \}_{p(k)}
$ and $b_j$ we get that elements $s$ and $t$ in $F^{\times }$
exist so that \par $ \{ b_1,...,b_j,...,b_k \}_{p(k)} = s \{
b_1,...,b_{j-1},b_{j+1},,...,b_k \}_{p(k-1)}b_j$
\par $ = st ( \{
b_1,...,b_{j-1},b_{j+1},,...,b_k \}_{p(k-1)}c_j)a_n$, \\ where $p(k-1)$ is a
corresponding vector prescribing an order of multiplications. \\ Again
applying the induction hypothesis to the product of $n-1$ elements $\{
b_1,...,b_{j-1},b_{j+1},,...,b_k \}_{p(k-1)}c_j$ we deduce that
there exists $w\in F^{\times }$ such that \par $ \{
a_{v(1)},...,a_{v(n)} \}_{u(n)} =
stw((...((a_1a_2)a_3)...)a_{n-1})a_n $.
\par Therefore a case remains when $m=1$. Let the first
multiplication in $ \{ a_{v(1)},...,a_{v(n)} \}_{u(n)}$ containing
$a_n$ be $(a_{v(k)}a_{v(k+1)})=: b_k$, consequently,
\par $ \{ a_{v(1)},...,a_{v(n)} \}_{u(n)} =  \{
b_{y(1)},...,b_{y(n-1)} \}_{w(n-1)} $ \\ for some $y\in S_{n-1}$ and
a vector $w(n-1)$ indicating an order of pairwise products, where
$b_j=a_{v(j)}$ for each $1\le j\le k-1$, also $b_{j-1}=a_{v(j)}$ for
each $k+1< j \le n$, where either $a_n=a_{v(k)}$ or
$a_n=a_{v(k+1)}$. From the induction hypothesis we deduce that there
exists $t_{n-1}\in F^{\times }$ so that
\par $t_{n-1}  \{ b_{y(1)},...,b_{y(n-1)} \}_{w(n-1)} = pb_k$ with
$ p = \{ b_1,...,b_{k-1}, b_{k+1},..., b_{n-1} \}_{w(n-1)} $.
Applying the induction hypothesis for $n=3$ we infer that there
exists $t_3\in F^{\times }$ such that
\par $ t_{n-1}t_3 \{ b_{y(1)},...,b_{y(n-1)} \}_{w(n-1)} =  (
pa)a_n$, \\ where either $a=a_{v(k+1)}$ or $a=a_{v(k)}$
correspondingly. From the induction hypothesis for $n-1$ it follows
that there exists $\tilde{t}_{n-1}\in F^{\times }$ so that
$\tilde{t}_{n-1} pa=(...((a_1a_2)a_3)...)a_{n-1}$ and hence $\{
a_1,...,a_n \}_{l(n)}= t_n \{ a_{v(1)},...,a_{v(n)} \}_{u(n)}$, \\
where $t_n =\tilde{t}_{n-1} t_{n-1}t_3 $.
\par $(ii)$. For a metagroup $G$ and the neutral element $v=id$
of the symmetric group $S_n$
we deduce analogously property $(1)$, since condition 1$(10)$ is
already unnecessary for $v=id$, where $id(k)=k$ for each $k\in {\bf
N}$.
\par {\bf 3. Lemma.} {\it If $G$ is a metagroup, then
for each $b\in G$ the identity is fulfilled $b\setminus e={\sf
t}_3((e/b),b,(b\setminus e))(e/b)$.}
\par {\bf Proof.} From formulas 1$(1)$, 1$(2)$ and 1$(9)$
we deduce that \par $e/b = (e/b)(b(b\setminus e)) = {\sf
t}_3^{-1}((e/b),b,(b\setminus e))((e/b)b)(b\setminus e)$\par $= {\sf
t}_3^{-1}((e/b),b,(b\setminus e))(b\setminus e)$.

\par {\bf 4. Definition.} For a metagroup $G$
and an associative unital ring $\cal T$
satisfying conditions $(1-3)$: \par $(1)$ $F^{\times }$ is contained in a center ${\cal C}({\cal T})$
of $\cal T$  \par $(2)$ $s_ja_j=a_js_j$ for each $j=1,...,n$, where $n$ is an arbitrary
natural number,
\par $(3)$ $s(ra)=(sr)a$ for each $s$ and $r$ in $\cal T$, and
$a\in G$, \\ by ${\cal T}[G]$ is denoted a
metagroup algebra over $\cal T$ of all formal sums $s_1a_1+...+s_na_n$,
where $s_1$,...,$s_n$ are in $\cal T$ and $a_1$,...,$a_n$ belong to
$G$, where \par $ ~ {\cal C}({\cal T}) := \{ b\in {\cal T}: \forall a\in {\cal T} ~ ab=ba \} $.

\par {\bf 5. Note.} Let $M$ be an additive commutative group such that
$M$ is a two-sided $G$-module, that is to each $g\in G$ there correspond
automorphisms $p(g)$ and $s(g)$ of $M$ putting for short $gx=p(g)x$
and $xg=xs(g)$ for each $g\in G$.
\par  Evidently, $M$ is a two-sided $G$-module if and only if it
is a two-sided ${\bf Z}[G]$-module according to the formulas
\par $(\sum_{g\in G}n(g)g)x=\sum_{g\in G}n(g)(gx)$ and
\par $x(\sum_{g\in G}n(g)g)=\sum_{g\in G}(xg)n(g)$,
\\ where $n(g)\in \bf Z$ for each $g\in G$, $ ~ {\bf Z}$ denotes the ring of all integers.
\par One can consider the additive group of integers $\bf Z$ as the trivial two-sided
$G$-module putting $gn=ng=n$ for each $g\in G$ and $n\in \bf Z$,
where $G$ is a metagroup.

\par {\bf 6. Examples.} Recall the following. Let $A$ be
 a unital algebra over a commutative associative unital ring $F$
supplied with a scalar involution $a\mapsto \bar{a}$ so that its
norm $N$ and trace $T$ maps have values in $F$ and fulfil
conditions:
\par $(1)$ $a\bar{a} =N(a)1$ with $N(a)\in F$,
\par $(2)$ $a+\bar{a}=T(a)1$ with $T(a)\in F$,
\par $(3)$ $T(ab)=T(ba)$
\\ for each $a$ and $b$ in $A$.
\par  If a scalar $f\in F$ satisfies the condition: $\forall a\in A ~ fa=0
\Rightarrow a=0$, then such element $f$ is called cancelable. For a
cancelable scalar $f$ the Cayley-Dickson doubling procedure provides
new algebra $C(A,f)$ over $F$ such that:
\par $(4)$ $C(A,f)=A\oplus Al$,
\par $(5)$ $(a+bl)(c+dl)=(ac - f\bar{d}b)+(da+b\bar{c})l$ and
\par $(6)$ $\overline{(a+bl)} = \bar{a} - bl$
\\ for each $a$ and $b$ in $A$. Then $l$ is called a doubling
generator. From definitions of $T$ and $N$ it follows that $\forall
a\in A, \forall b\in A$ $~T(a)=T(a+bl)$ and $N(a+bl)=N(a)+ fN(b)$.
The algebra $A$ is embedded into $C(A,f)$ as $A\ni a\mapsto (a,0)$,
where $(a,b)=a+bl$. Put by induction $A_n(f_{(n)}) =
C(A_{n-1},f_n)$, where $A_0=A$, $f_1=f$, $n=1, 2,...$, $ ~
f_{(n)}=(f_1,...,f_n)$. Then $A_n(f_{(n)})$ are generalized
Cayley-Dickson algebras, when $F$ is not a field, or Cayley-Dickson
algebras, when $F$ is a field.
\par It is natural to put
$A_{\infty }(f) := \bigcup_{n=1}^{\infty }A_n(f_{(n)})$, where
$f=(f_n: n\in {\bf N})$.  If $char (F)\ne 2$, let $Im (z)=z- T(z)/2$
be the imaginary part of a Cayley-Dickson number $z$ and hence
$N(a):=N_2(a,\bar{a})/2$, where $N_2(a,b):=T(a\bar{b})$.
\par If the doubling procedure starts from $A=F1=:A_0$, then
$A_1=C(A,f_1)$ is a $*$-extension of $F$. If $A_1$ has a basis $ \{
1, u \} $ over $F$ with the multiplication table $u^2=u+w$, where
$w\in F$ and $4w+1\ne 0$, with the involution $\bar{1}=1$,
$\bar{u}=1-u$, then $A_2$ is the generalized quaternion algebra,
$A_3$ is the generalized octonion (Cayley-Dickson) algebra.
\par When $F=\bf R$ and $f_n=1$ for each $n$ by ${\cal A}_r$ will be denoted
the real Cayley-Dickson algebra with generators $i_0,...,i_{2^r-1}$
such that $i_0=1$, $~i_j^2=-1$ for each $j\ge 1$, $~i_ji_k=-i_ki_j$
for each $j\ne k\ge 1$.  Frequently $\bar{a}$ is also denoted by
$a^*$ or $\tilde a$. \par Let $A_n$ be a Cayley-Dickson algebra over
a commutative associative unital ring $\cal R$ of characteristic
different from two such that $A_0=\cal R$, $n\ge 2$. Take its basic
generators $i_0, i_1,..., i_{2^n-1}$, where $i_0=1$. Choose
$F^{\times }$ as a multiplicative subgroup contained in the ring
$\cal R$ such that $f_j\in F^{\times }$ for each $j=0,...,n$. Put
$G_n = \{ i_0, i_1,...,i_{2^n-1} \} \times F^{\times }$. Then $G_n$
is a central metagroup. \par More generally let $H$ be a group such
that $F^{\times }\subset H$, with relations $hi_k=i_kh$ and
$(hg)i_k=h(gi_k)$ for each $k=0,1,...,2^n-1$ and each $h$ and $g$ in
$H$. Then $G_n = \{ i_0, i_1,...,i_{2^n-1} \} \times H$ is also a
metagroup. The latter metagroup is noncentral, when $H$ is
noncommutative. Analogously to the Cayley-Dickson algebra $A_{\infty
}(f)$ a metagroup $G_{\infty }$ corresponds, when $n=\infty $.
\par Generally metagroups need not be central. From given
metagroups new metagroups can be constructed using their direct or
semidirect products. Certainly each group is a metagroup also.
Therefore, there are abundant families of noncentral metagroups and also of central metagroups
different from groups.
\par In another way smashed products of groups can be considered analogously to the Cayley-Dickson
doubling procedure providing another examples of metagroups.

\par {\bf 7. Definitions.} If $\cal R$ is a ring, which may be nonassociative relative to the
multiplication, a two-sided module $M$ over $\cal R$ is called cyclic, if an element
$y\in M$ exists such that $M={\cal R}(y{\cal R}) = \{ s(yp): ~ s, p\in {\cal R} \} $
and \par $M=({\cal R}y){\cal R} = \{ (sy)p: ~ s, p\in {\cal R} \} $.
\par Take an algebra $A={\cal T}[G]$ and a two-sided $A$-module $M$, where
$\cal T$ is an associative unital ring satisfying conditions 4$(1-3)$.
Let $M$ have the decomposition $M=\sum_{g\in G}M_g$ as a two-sided ${\cal T}$-module,
where $M_g$ is a two-sided ${\cal T}$-module for each $g\in G$, $G$ is a metagroup,
and let $M$ satisfy the following conditions:
\par $(1)$ $hM_g=M_{hg}$ and $M_gh=M_{gh}$,
\par $(2)$ $(bh)x_g=b(hx_g)$ and $x_g(bh)=(x_gh)b$ and $bx_g=x_gb$,
\par $(3)$ $(hs)x_g={\sf t}_3(h,s,g) h(sx_g)$ and $(hx_g)s= {\sf t}_3(h,g,s) h(x_gs)$
and \par $(x_gh)s={\sf t}_3(g,h,s) x_g(hs)$ \\ for every $h, g, s$ in $G$ and $b\in \cal T$
and $x_g\in M_g$. Then a two-sided $A$-module $M$ will be called $G$-graded.

\par {\bf 8. Example.} Let $\cal T$ be a commutative associative unital ring,
let also $G$ be a metagroup and $A={\cal T}[G]$ be a metagroup algebra (see \S 4), where
$A$ is considered as a $\cal T$-algebra.
Put $K_{-1}= A$, $K_0 = A \otimes_{\cal T} A$ and by induction
$K_{n+1} = K_n\otimes_{\cal T} A$ for each natural number $n$.
Each $K_n$ is supplied with a two-sided $A$-module structure:
\par $(1)$ $p\cdot (x_0,...,x_{n+1})= ((px_0),...,x_{n+1})=...
=$\par $(x_0,...,(px_{n+1}))= (x_0,...,x_{n+1})\cdot p$ \par for each $p\in {\cal T}$,
where $0\cdot (x_1,...,x_n)=0$,
\par $(2)$ $(xy)\cdot (x_0,...,x_{n+1})={\sf t}_3\cdot (x\cdot (y\cdot (x_0,...,x_{n+1})))$ \par with
${\sf t}_3={\sf t}_3(x,y,b)$, (see also formula 1$(9)$ above);
\par $(3)$ ${\sf t}_3\cdot ((x_0,...,x_{n+1})\cdot (xy)) = ((x_0,...,x_{n+1})\cdot x)\cdot y$
with ${\sf t}_3={\sf t}_3(b,x,y)$;
\par $(4)$ $(x\cdot (x_0,...,x_{n+1}))\cdot y = {\sf t}_3 \cdot (x\cdot ((x_0,...,x_{n+1})\cdot y))$ with
${\sf t}_3={\sf t}_3(x,b,y)$
\par $(5)$ $x\cdot (x_0,...,x_{n+1})= t_{n+3} (x,x_0,...,x_{n+1};v_0(n+3);l(n+3))\cdot ((xx_0),x_1,...,x_{n+1})$
\\ where $\{ x, x_0,...,x_{n+1} \} _{v_0(n+3)} = x \{ x_0,...,x_{n+1} \} _{l(n+2)}$,
\par $\{ x_0,...,x_{n+1} \}_{l(n+2)} = \{ x_0,...,x_n \} _{l(n+1)} x_{n+1}$, \par $\{ x_0 \}_{l(1)}= x_0$,
$ \{ x_0 x_1 \} _{l(2)} =x_0x_1$; \\ where $b= \{ x_0,...,x_{n+1} \} _{l(n+2)}$,
\par $t_{n} (x_1,...,x_{n};u(n),w(n)) := t_{n} (x_1,...,x_{n};u(n), w(n)|id)$ \\ using shortened notation;
\par $(6)$ $(x_0,...,x_{n+1})\cdot x= t_{n+3}(x_0,...,x_{n+1},x;l(n+3),v_{n+2}(n+3))\cdot (x_0,...,x_n, (x_{n+1}x))$
\\ for every $x, y, x_0,...,x_{n+1}$ in $G$, where $(x_0,...,x_{n+1})$ denotes a
basic element of $K_n$ over $\cal T$ corresponding to the left ordered tensor product
\par $(...((x_0\otimes x_1)\otimes x_2)...\otimes x_n)\otimes x_{n+1}$,
\par $\{ x_0,...,x_{n+1},x \} _{v_{n+2}(n+3)} = \{ x_0,...,x_n,x_{n+1}x \} _{l(n+2)}$.

\par {\bf 9. Proposition.} {\it For each metagroup algebra $A={\cal T}[G]$
(see \S 4) an acyclic left $A$-complex $\cal K$ exists.}
\par {\bf Proof.} Take two-sided $A$-modules $K_n$ as in example 8.
We define a boundary $\cal T$-linear operator
$\partial _n: K_n\to K_{n-1}$ on $K_n$. On
basic elements we put it to be given by the formulas:
\par $(1)$ $\partial _n((x\cdot (x_0,x_1,...,x_n,x_{n+1}))\cdot y) =$\par $\sum_{j=0}^{n}
(-1)^{j}\cdot t_{n+4}(x,x_0,...,x_{n+1},y;l(n+4),u_{j+1}(n+4))$\par $\cdot
((x\cdot(<x_0, x_1,..., x_{n+1}
>_{j+1,n+2}))\cdot y) $, where
\par $(2)$ $ <x_0,...,x_{n+1} >_{1,n+2} := ((x_0x_1),x_2,...,x_{n+1})$,
\par $(3)$ $ <x_0,...,x_{n+1} >_{2,n+2} :=  (
x _0,(x_1x_2),x_3,...,x_{n+1} )$,...,
\par $(4)$ $ < x_0,...,x_{n+1} >_{n+1,n+2} :=  (
x_0,...,x_{n-1},(x_nx_{n+1}) )$,
\par $(5)$ $\partial _0 (x\cdot (x_0,x_1))\cdot y =(x\cdot (x_0x_1))\cdot y$,
\par $(6)$ $\{ x_0,x_1,...,x_{n+1} \}_{l(n+2)}
:= (...((x_0x_1)x_2)...)x_{n+1} $;
\par $(7)$ $ \{ x, x_0, ...,x_{n+1}, y \}_{u_1(n+4)} :=  (x \{ (x_0x_1),
x_2,...,x_{n+1} \}_{l(n+1)})y$,...,
\par $(8)$ $ \{ x, x_0,...,x_{n+1},y \}_{u_{n+1}(n+4)} :=  (x \{
x_0,x_1,...,(x_nx_{n+1}) \}_{l(n+1)})y$
\\ for each $x, x_0,...,x_{n+1}, y$ in $G$.
On the other hand, from formulas 1$(1)$ and 1$(2)$ it follows that
$t_{n+4}(x,x_0,...,x_{n+1},y;l(n+4),u_{j+1}(n+4)) =$\par $ t_{n+2}
(x_0,...,x_{n+1};l(n+2),v_{j+1}(n+2))$  for each $j=0,...,n$, where
\par $(9)$ $ \{ x_0, ...,x_{n+1} \}_{v_1(n+2)} :=  \{ (x_0x_1),
x_2,...,x_{n+1} \}_{l(n+1)}$,...,
\par $(10)$ $ \{ x_0,...,x_{n+1} \}_{v_{n+1}(n+2)} :=  \{
x_0,x_1,...,(x_nx_{n+1}) \}_{l(n+1)}$
\\ for every $x_0,...,x_{n+1}$ in $G$.
Therefore, $\partial _n$ is a left and right $A$-homomorphism of $(A,A)$-modules.
In particular, \par $\partial
_1((x\cdot (x_0,x_1,x_2))\cdot y) = (x\cdot ((x_0 x_1),x_2))\cdot y - {\sf t}_{3}(x_0,x_1,x_2)
\cdot (x\cdot (x_0, (x_1x_2)))\cdot y $,
\par $\partial _2((x\cdot (x_0,x_1,x_2,x_3))\cdot y)= (x\cdot ((x_0x_1),x_2,x_3))\cdot y
-t_4(x_0,...,x_3; l(4),v_2(4))\cdot ((x\cdot
(x_0,(x_1x_2),x_3))\cdot y) + t_4(x_0,...,x_3; l(4),v_3(4);id)
((x\cdot (x_0,x_1,(x_2x_3)))\cdot y)$.
\par Define a $\cal T$-linear homomorphism
$s_n: K_n \to K_{n+1}$, which on basic elements has the form:
\par $(11)$ $s_n(x_0,...,x_{n+1})= (1,x_0,....,x_{n+1})$
\\ for every $x_0,...,x_{n+1}$ in $G$.
From formulas 1$(9)$, 1$(10)$ and 2$(1)$ the identities
\par $(12)$ $t_n(x_1,...,x_n;q(n),u(n)|v)
t_n(x_1,...,x_n;u(n),q(n)|v^{-1})=1$
\par $(13)$ $t_n(x_1,...,x_n;q(n),u(n)) t_n(x_1,...,x_n;u(n),w(n))$ \par $=
t_n(x_1,...,x_n;q(n),w(n))$
\\ follow for every elements
$x_1$,...,$x_n$ in the metagroup $G$, vectors $q(n)$, $u(n)$ and $w(n)$
indicating orders of their multiplications, $v\in S_n$ and $n\in \bf
N$. Moreover, \par $(14)$
$t_{n+1}(1,x_1,...,x_n;q(n+1),u(n+1)|v(n+1))$ \par $=
t_n(x_1,...,x_n;q(n),u(n)|v(n))$ \\ for $q(n)$, $u(n)$ and $v(n)$
obtained from $q(n+1)$, $u(n+1)$ and $v(n+1)$ correspondingly by
taking into account the identity $1b=b1=b$ for each $b\in G$.
Hence \par $s_n((x_0,...,x_{n+1})\cdot y) = (s_n(x_0,...,x_{n+1}))\cdot y$
\\ for every $x_0,...,x_{n+1}, y$ in $G$.
\par Let $p_n: K_{n+1}\to K_n$ be a $\cal T$-linear mapping such that
\par $(15)$ $p_n(a\otimes b)=a\cdot b$ and $p_n(b\otimes a)=b\cdot a$ for each
$a\in K_n$ and $b\in A$. Therefore, from formulas $(13)$ and $(14)$ we deduce that
$p_ns_n=id$ is the identity on $K_n$, consequently,
$s_n$ is a monomorphism.
\par Therefore, from formulas 4$(1-3)$, 7$(1-3)$, $(1)$ and $(9,11,13,14)$ we infer
that
\par $(\partial _{n+1} s_n + s_{n-1} \partial _n)(x_0,...,x_{n+1}) =$
\par $= \partial _{n+1} (1,x_0,...,x_{n+1}) +$ \par $s_{n-1}(\sum_{j=0}^{n}
(-1)^{j} t_{n+2}(x_0,...,x_{n+1};l(n+2),v_{j+1}(n+2))\cdot $\par $ <x_0,
..., x_{n+1} >_{j+1,n+2}))=$
\par $=\sum_{j=0}^{n+1}
(-1)^{j}t_{n+3}(1,x_0,...,x_{n+1};l(n+3),v_{j+1}(n+3))\cdot $\par $
<1,x_0, x_1,..., x_{n+1} >_{j+1,n+3}+$
\par
$+ \sum_{j=0}^{n}
(-1)^{j} t_{n+2}(x_0,...,x_{n+1};l(n+2),v_{j+1}(n+2))\cdot $\par $
<1,x_0, ...,x_{n+1} >_{j+2,n+3} = (x_0,...,x_{n+1})$,\\
for every $x_0$,...,$x_{n+1}$ in $G$.
\par Thus the homotopy conditions
\par $(16)$ $\partial _{n+1}s_n+s_{n-1}\partial _n =1$ for each $n\ge 0$ \\
are fulfilled, where $1$ denotes the identity operator on $K_n$.
Therefore the recurrence relation
\par $\partial _n\partial _{n+1} s_n = \partial _n
(1-s_{n-1}\partial _n) = \partial _n - (\partial _n s_{n-1})\partial
_n = $\par $ \partial _n - (1-s_{n-2}\partial _{n-1})\partial _n = s_{n-2}
\partial _{n-1} \partial _n$
\\ is accomplished. On the other hand, from Formula $(11)$ it follows that $K_{n+1}$ as the left $A$-module is
generated by $s_nK_n$. Then proceeding by induction in $n$ we deduce that
$\partial _n \partial _{n+1}=0$ for each $n\ge 0$, since $\partial _0 \partial _1=0$
according to formulas $(1)$ and $(5)$. Mention also that $K_0=A\otimes _{\cal T}A$ coincides
with $A\otimes _{\cal T}A^{op}$ as
a left and right $A$-module, hence the mapping $\partial _0: K_0\to K_{-1}$
provides the augmentation $\epsilon : A^e\to A$, where $A^e := A\otimes _{\cal T}A^{op}$
denotes the enveloping algebra of $A$, $A^{op}$ notates the opposite algebra of $A$.
Thus identities $(16)$ mean that the left complex $\cal K$
\par  $0\leftarrow A \mbox{ }_{
\overleftarrow{\partial _0}} K_0 \mbox{ }_{\overleftarrow{\partial _1}} K_1
\mbox{ }_{\overleftarrow{\partial _2}} K_2 \overleftarrow ...
\mbox{ }_{\overleftarrow{\partial _n}} K_n \mbox{ }_{\overleftarrow{\partial _{n+1}}} K_{n+1} \leftarrow ...$.
is acyclic.

\par {\bf 10. Examples. 1.} For the Cayley-Dickson algebra $A_n$ over a field $F$
of characteristic not equal to two let $G=G_n$ as the
(multiplicative) metagroup consist of all elements $bi_k$ with $b\in
F^{\times }$, $k=0,1,2,...$, where $ i_0, i_1, i_2,...$ are
generators of the Cayley-Dickson algebra $A_n$, $2\le n\le \infty $.
Then $M=A_n^j$ is the module over ${\bf Z}[G]$, where $j\in {\bf
N}$.
\par {\bf 10.2.} For a topological space $U$ it is possible to consider the
module $M=C(U,A_n^j)$ of all continuous mappings from $U$ into
$A_n^j$, $j\in \bf N$, $A_n^j$ is supplied with the box product
topology.
\par {\bf 10.3.} If $(U,{\cal B}, \mu )$ is a measure space, where $\mu :
{\cal B}\to [0, \infty )$ is a $\sigma $-additive measure on a
$\sigma $-algebra ${\cal B}$ of a set $U$, for ${\bf F}={\bf R}$ and
$f_k=1$ for each $k$, it is possible to consider the space
$L_p((U,{\cal B}, \mu ),A_n^j)$ of all $L_p$ mappings from $U$ into
$A_n^j$, where $A_n$ is taken relative to its norm induced by the
scalar product $Re (\bar{y}z) =(y,z)$, $j\in \bf N$, $1\le p\le
\infty $.
\par {\bf 10.4.} For an additive group $H$ one can consider the trivial action of $A$
on $H$. Therefore, the direct product $M\bigotimes H$ becomes an
$A$-module for an $A$-module $M$. In particular, $H$ may be a ring.
\par {\bf 10.5.} If there is another ring $\cal S$ and a homomorphism
$\phi : {\cal S}\to {\cal T}$, then each left (or right) $\cal T$-module $M$
can be considered as a left (or right correspondingly) $\cal S$-module by the rule:
\par $bm=(\phi b)m$ (or $mb=m(\phi b)$ correspondingly) for each $b\in {\cal S}$ and $m\in M$.
\par Vice versa if $M$ is a right (or left) ${\cal S}$-module,
then there exist the right (or left correspondingly) module
$M_{(\phi )} = M\otimes _{\cal S}{\cal T}$ (or $\mbox{ }_{(\phi )}M
= {\cal T}\otimes _{\cal S} M$ correspondingly) called the right (or
left correspondingly) covariant $\phi $-extension of $M$. Similarly
are defined the contravariant right and left extensions $M^{(\phi
)}=Hom_{\cal S}({\cal T},M)$ or $\mbox{ }^{(\phi )}M$ for a right or
left $\cal S$-module $M$ respectively.
\par This also can be applied to a metagroup algebra
$A={\cal S}[G]$ over a commutative associative unital ring $\cal S$ as in Examples 6.
Then changing a ring we get right $A_{(\phi )} $ or $A^{(\phi )} $
and left $\mbox{ }_{(\phi )}A$ or $\mbox{ }^{(\phi )}A$
algebras over $\cal T$. Then imposing
the relation $ta=at$ for each $a\in A$ and $t\in T$ provides
a metagroup algebra over $\cal T$ which also has a two-sided ${\cal T}$-module
structure. It will be denoted by $\mbox{ }_{(\phi )}A_{(\phi )}$ or $\mbox{ }^{(\phi )}A^{(\phi )}$
respectively.
Particularly, this is applicable to cases when ${\bf Z}[F^{\times }]
\subset {\cal S}$ or $\phi $ is an embedding.

\par {\bf 11. Notation.} Let $A={\cal T}[G]$ be a metagroup algebra (see \S 4).
Put $L_0={\cal T}$, $L_1=A$ and by induction $L_{n+1}=L_n\otimes
_{\cal T}A$ for each natural number $n$. \par  If $N$ is a two-sided
$A$-module it can also be considered as a left $A^e$-module by the
rule: $(x\otimes y^*)b:= (x\otimes b)\otimes y$ for each $x\in A$,
$y^*\in A^{op}$ and $b\in N$, $ ~ A^e=A\otimes _{\cal T}A^{op}$ is
an enveloping algebra, where $A^{op}$ denotes the opposite algebra
of $A$, where $y^*$ in $A^{op}$ is the corresponding to $y$ in $A$
element.

\par {\bf 12. Proposition.} {\it If $\cal K$ is an acyclic left $A$-complex
for a metagroup algebra $A={\cal T}[G]$ as in Proposition 9 and $M$
is a two-sided $A$-module satisfying conditions 7$(1-3)$, then there exists a cochain
complex $Hom ({\cal L},M)$:
\par $(1)$ $0\to Hom_{\cal T}(L_0,M)_{\overrightarrow{{\epsilon
^*}}}Hom_{\cal T}(L_1,M)_{\overrightarrow{{\delta
^1}}}Hom_{\cal T}(L_2,M)_{\overrightarrow{{\delta ^2}}}$\par
$Hom_{\cal T}(L_3,M)_{\overrightarrow{{\delta
^3}}}Hom_{\cal T}(L_4,M)_{\overrightarrow{{\delta ^4}}}...$ \\ such that
$f\in Hom_{\cal T}(L_1,M)$ is a cocycle if and only if $f$ is a $\cal T$-linear
derivation from $A$ into $M$.}
\par {\bf Proof.} The notations of sections 8 and 11 permit to write
each basic element $(x_0,...,x_{n+1})$ of $K_n$ over ${\cal T}$ as \par $(1)$ $(x_0,...,x_{n+1})=
t_{n+2}(x_0,...,x_{n+1};l(n+2),w(n+2))\cdot $
\par $((x_0\otimes (x_1,...,x_n))\otimes x_{n+1})$ and \par $(2)$
$(x_0,...,x_{n+1})=
t_{n+2}(x_0,...,x_{n+1};l(n+2),w(n+2))\cdot $\par $
(z\otimes (x_1,...,x_n))$, \\ where $(x_1,...,x_n)$ is a basic element
in $L_n$ for  every $x_0,...,x_{n+1}$ in $G$, \par $ \{ x_0,...,x_{n+1} \} _{w(n+2)}=
(x_0 \{ x_1,...,x_n \} _{l(n)})x_{n+1}$, \\ $z\in A^e$, $z=x_0\otimes x^*_{n+1}$.
\par Each homomorphism $f\in Hom_{\cal T}(L_n,M)$ is characterized by its
values on elements $(x_1,...,x_n)$, where $x_1$,..,$x_n$ belong to a metagroup
$G$. Consider $f$ as a $\cal T$-linear function from $A^n$ into $M$. Since
$M$ satisfies conditions 7$(1-3)$, then $f$ has the decomposition
\par $(2)$ $f(x_1,...,x_n)=\sum_{g\in G} f_g(x_1,...,x_n)$,\\ where
$f_g: G^n \to M_g$ for every $g$ and $x_1,...,x_n$ in $G$.
\par Therefore, the restrictions follow from conditions 7$(1-3)$, which
take into account the nonassociativity of $G$:
\par $(3)$ $(xy)\cdot f_g(x_1,...,x_n)={\sf t}_3(x,y,g)\cdot (x\cdot (y\cdot
f_g(x_1,...,x_n)))$,
\par $(4)$ ${\sf t}_3(g,x,y)\cdot (f_g(x_1,...,x_n)\cdot (xy)) = (f_g(x_1,...,x_n)\cdot x)\cdot
y$,
\par $(5)$ $(x\cdot f_g(x_1,...,x_n))\cdot y = {\sf t}_3(x,g,y) \cdot (x\cdot (f_g(x_1,...,x_n)\cdot y))$
\\ for every $g$ and $x, y, x_1,..., x_n$ in $G$, where coefficients ${\sf t}_3$ are prescribed
by formula 1$(9)$, also \par $(6)$ $x\cdot f_g(x_1,...,x_n):=x\cdot
(f_g(x_1,...,x_n))$ and \par $(7)$ $f_g(x_1,...,x_n)\cdot y :=
(f_g(x_1,...,x_n))\cdot y$.
\par For $n=0$ and $g=e$
naturally the identities are fulfilled:
\par $(8)$ $(xy)\cdot f_e( ~ )=x\cdot (y\cdot f_e( ~ ))$, $(f_e( ~ )\cdot
x)\cdot y=f_e( ~ )\cdot (xy)$ and $(x\cdot f_e( ~ ))\cdot y=x\cdot (f_e( ~
)\cdot y)$.

\par  We define a coboundary operator taking into account
the nonassociativity of the (multiplicative) metagroup $G$:
\par $(9)$ $(\delta ^nf)(x_1,...,x_{n+1}) = $\par $\sum_{j=0}^{n+1} (-1)^{j}
t_{n+1}(x_1,...,x_{n+1};l(n+1),u_{j+1}(n+1))\cdot [f,x_1, x_2,...,
x_{n+1} ]_{j+1,n+1}$,  where
\par $(10)$ $ [f,x_1,...,x_{n+1}]_{1,n+1} :=
x_1 \cdot f(x_2,...,x_{n+1} )$,
\par $\{ x_1,...,x_{n+1} \}_{u_1(n+1)} = x_1 \{ x_2,...,x_{n+1} \}_{l(n)}$;
\par $(11)$ $ [f,x_1,...,x_{n+1}]_{2,n+1} :=
f((x_1x_2),...,x_{n+1} )$,
\par $\{ x_1,...,x_{n+1} \}_{u_2(n+1)} = \{ (x_1x_2),...,x_{n+1} \}_{l(n)}$;...;
\par $(12)$ $ [f,x_1,...,x_{n+1}]_{n+1,n+1} :=  f(
x_1,x_2,...,(x_nx_{n+1}) )$;
\par $\{ x_1,...,x_{n+1} \}_{u_{n+1}(n+1)} = \{ x_1,...,(x_nx_{n+1}) \}_{l(n)}$;
\par $(13)$ $ [f,x_1,...,x_{n+1}]_{n+2,n+1} :=  f(
x_1,x_2,...,x_n)\cdot x_{n+1}$,
\par $\{ x_1,...,x_{n+1} \}_{u_{n+2}(n+1)} =  \{ x_1,...,x_{n+1} \}_{l(n+1)}
=(...((x_1x_2)x_3)...x_n)x_{n+1}$;
\\ with $u_0(n+1)=l(n+1)$.
From $G^{n+1}$ onto $K_{n+1}$ the homomorphism $(\delta ^nf)$ is
extended by $\cal T$-linearity. On the other hand, condition 7$(1)$
implies that
\par $(14)$ for each $b\in G$ there exists $h_{1,b}$ so that
$h_{1,b}: K_{n+1}\to M_1$ and $f_b=h_{1,b}L_b$, where
$L_b$ is the left multiplication operator on $b$:
\par $(15)$
$(h_{1,b}L_b)(x_1,...,x_n)= b\cdot (h_{1,b}(x_1,...,x_n))$\\
for every $x_1,...,x_n$ in $G$.
Moreover, $zg=0$ (or $gz=0$) in ${\bf Z}[G]$ for $g\in G$ and $z\in
{\bf Z}[G]$ if and only if $z=0$, since $G$ is a metagroup. \par In
virtue of Proposition 9 these formulas imply that $\delta
^{n+1}\circ \delta ^n=0$ for each $n$, since \par $(\delta ^{n+1}\circ
\delta ^nf)(x_1,...,x_{n+2})=f(\partial _{n-1}\circ
\partial _{n}(x_1,...,x_{n+2}))$ \\ for every $x_1,...,x_{n+2}$ in
$G$. Thus the complex given by formula $(1)$ is exact.
\par Particularly, $f\in Hom_{\cal T}(L_0,M)$ is a cocycle if and only if
\par $(16)$ $(\delta ^0f)(x)=xf(~)-f(~)x=0$ for each $x\in G$.
\\ Mention that $Hom _{\cal T}(L_0,M)$ is isomorphic with $M$.
\par One dimensional cochain $f\in Hom_{\cal T}(L_1,M)$ is determined by
a mapping $f: G\to M$. Taking into account formula
$(9)$ we infer that it is a cocycle if and only if
\par $(17)$ $t_{2}(x,y;l(2),u_{1}(2))\cdot x\cdot f(y)-$\par $
t_{2}(x,y;l(2),u_{2}(2))\cdot
f(xy)+t_{2}(x,y;l(2),u_{3}(2))\cdot f(x)\cdot y $\par $ =x\cdot
f(y)-f(xy)+f(x)\cdot y =0$
\\ for each $x$ and $y$ in $G$. That is $f$ is a derivation
from the metagroup $G$ into the $G$-module $M$. There is the
embedding ${\cal T}\hookrightarrow A$ as ${\cal T}e$, since $e=1\in G$.
Thus $f$ has a $\cal T$-linear extension to a $\cal T$-linear
derivation from $A$ into $M$:
\par $(18)$ $f(xy)=x\cdot f(y)+f(x)\cdot y$.

\par {\bf 13. Remark.} Suppose that conditions of Proposition 12 are fulfilled.
A two-dimensional cochain is a $2$-cocycle, if and only if
\par $(\delta ^2f)(x_1,x_2,x_3) = $\par $\sum_{j=0}^{3} (-1)^{j}
t_{3}(x_1,x_2,x_3;l(3),u_{j+1}(3))\cdot [f,x_1, x_2,x_3]_{j+1,3}$
\par $=  t_{3}(x_1,x_2,x_3;l(3),u_{1}(3))\cdot x_1\cdot f(x_2,x_3) -$
\par $t_{3}(x_1,x_2,x_3;l(3),u_{2}(3))\cdot f((x_1x_2),x_3) +$
\par $t_{3}(x_1,x_2,x_3;l(3),u_{3}(3))\cdot f(x_1,(x_2x_3)) -$
\par $t_{3}(x_1,x_2,x_3;l(3),u_{4}(3))\cdot f(x_1,x_2)\cdot
x_3=0$ \par that is
\par $(1)$ ${\sf t}_{3}(x_1,x_2,x_3)\cdot x_1\cdot
f(x_2,x_3)+ {\sf t}_3(x_1,x_2,x_3)\cdot f(x_1,(x_2x_3))
 $\par $ =f((x_1x_2),x_3)+f(x_1,x_2)\cdot x_3$ \\ for each $x_1, x_2$ and
$x_3$ in $G$. \par As usually $Z^n(A,M)$ denotes the set of all
$n$-cocycles, the notation $B^n(A,M)$ is used for the set of
$n$-coboundaries in $H_{\cal T}(L_n,M)$. Since as the additive group $M$ is
commutative, then there are defined groups of cohomologies $H^n(A,M)
= Z^n(A,M)/B^n(A,M)$ as the quotient (additive) groups.
\par For $n=0$ coboundaries are put zero and hence $H^0(A,M) \cong M^A$.
In the case $n=1$ a mapping $f: A\to M$ is a coboundary if there
exists an element $m=h(~)\in M$ for which $f(x)=xm-mx$. Such
derivation $f$ is called an inner derivation of $A$ defined by an
element $m\in A$. The set of all inner derivations is denoted by
$Inn_{\cal T}(A,M)$. From the cohomological point of view the
additive group $H^1(A,M)$ is interpreted as the group of all outer
derivations $H^1(A,M)\cong Out_{\cal T} (A,M) \cong Der_{\cal T}
(A,M)/Inn_{\cal T} (A,M)$, where $Z^1(A,M)=Der_{\cal T}(A,M)$,
$Inn_{\cal T}(A,M)=B^1(A,M)$, where the family of all derivations
($\cal T$-homogeneous derivations) from $X$ into a two-sided module
$M$ over $\cal T$ is denoted by $Der(X,M)$ (or $Der_{\cal T} (X,M)$
respectively).
\par A two-cochain $f: G\times G\to M$ is a two-coboundary, if
an one-cochain $h: G\to M$ exists such that for each $x$ and $y$ in
$G$ the identity is fulfilled:
\par $(2)$ $f(x,y) = (\delta h)(x,y) =\sum_{j=0}^{2} (-1)^{j}
t_{2}(x,y;l(2),u_{j+1}(2))\cdot [h,x_1, x_2]_{j+1,2}$, \par $=
x\cdot h(y) - h(xy) +h(x)\cdot y$.

\par Let $A={\cal T}[G]$ be a metagroup algebra over a commutative associative unital ring
$\cal T$ (see also sections 1 and 4).
\par Let $M$, $N$ and $P$ be left $A$-modules and a short exact
sequence exists:
\par $(3)$ $0\to M_{\overrightarrow{{\xi }}}
P_{\overrightarrow{{\eta }}}N\to 0$,
\\ where $\xi $ is an embedding. Then $P$ is called an
enlargement of a left $A$-module $M$ with the help of
a left $A$-module $N$. If there is another enlargement of $M$ with the
help of $N$:
\par $(4)$ $0\to M_{\overrightarrow{{\xi '}}}
{P'}_{\overrightarrow{{\eta '}}}N\to 0$ \\ such that an isomorphism
$\pi : P\to P'$ exists for which $\pi \xi =\xi '1_M$ and $1_N\eta =
\eta ' \pi $, then enlargements $(3)$ and $(4)$ are called
equivalent, where $1_M: M\to M$ notates the identity mapping,
$1_M(m)=m$ for each $m\in M$.
\par In the particular case, when $P=M\oplus N$, also $\xi $ is an
identifying mapping with the first direct summand and $\eta $ is a
projection on the second direct summand, an enlargement is called
trivial.

\par {\bf 14. Theorem.} {\it Let a nonassociative algebra $A$ and left $A$-modules $M$ and $N$ be
as in section 13. Then the family $T=Hom_{\cal T}(N,M)$ can be supplied
with a two-sided $A$-module structure such that $H^1(A,T)$ is the set of
classes of modules $M$ with the quotient module $N$.}

\par {\bf Proof.} The family $T$ is a left module over a ring $\cal T$
and it can be supplied with a two-sided $A$-module structure:
\par $(1)$ $\forall ~ r\in T$ and $\forall ~ n\in N$ and $\forall ~ a\in A$:
\par $(a\cdot r)(n) = a\cdot (r(n))$ and $(r\cdot a)(n) =
r(a\cdot n)$.
\par In virtue of Proposition 12 each element $f\in Z^1(A,T)$
defines a (generalized) derivation by formula 12$(18)$. Each zero
dimensional cochain $m\in M$ provides an inner derivation
$\delta ^0m(a)=am-ma$ due to formula 12$(16)$. Then an one cocycle $f$ induces
an enlargement by formula 13$(3)$ with $P=M\oplus N$ being the
direct sum of left $A$-modules in which $N$ is a submodule and
with the left action of $A$ on $N$: $a\circ n = a\cdot n+f(a)\cdot
n$ for each $n\in N$ and $a\in A$. Suppose that a class of an one
cocycle $f$ is zero, that is an element $u\in T$ exists so that
$f=\delta ^1u$. Then elements of the form $m+u(m)$ form its
submodule $M'$ isomorphic with $M$. Moreover, $P=N\oplus M'$ is the
direct sum of $A$-modules. Thus an enlargement is trivial.
\par Vice versa. Suppose that an enlargement given by formula
13$(3)$ exists so  that a $\cal T$-linear mapping $\gamma : N\to P$
exists satisfying the restriction $\eta \gamma =1_N$, where $1_N$ notates the
identity mapping on $N$. Then for each
element $a$ of the algebra $A$ we put $f(a)n=\gamma (a\cdot n) -
(a\cdot \gamma )(n)$ for all $n\in N$.  Therefore $\gamma \in
Z^1(A,T)$. When $f$ is an inner derivation, a clefting follows
analogously to the proof above. On the other hand, when an
enlargement clefts, that is there exists an $A$-homomorphism $\omega
: N\to P$ fulfilling the restriction $\eta \omega =1_N$, then we
take $u\in T$ such that $u(n)=\omega (n)-\gamma (n)$.

\par {\bf 15. Theorem.} {\it Suppose that $A$ is a nonassociative metagroup algebra
over a commutative associative unital ring $\cal T$,
a left $A$-module $N$ and a two-sided $A$-module $M$ are given. Then
$H^2(A,T)$ is the set of classes of enlargements of $A$ with a
kernel $M$ such that $M^2 = \{ 0 \} $ and with the quotient algebra
$A$. Moreover, an action of $A$ on $M$ in this enlargement coincides
with the structure of a two-sided $A$-module on $M$.}
\par {\bf Proof.} If $P$ is an enlargement with a kernel $M$ such
that $M^2=0$ and a quotient module $A=P/M$ and $a=p+M$ with $p\in
P$, then $a\cdot m=p\cdot m$ and $m\cdot a=m\cdot p$ supply $M$ with
the two-sided $A$-module structure. Take a $\cal T$-linear mapping $\gamma
:A\to P$ inverse from the left to a natural epimorphism and put
$f(a,b)=\gamma (ab)-\gamma (a)\gamma (b)$ for each $a$ and $b$ in
$A$. Then we infer that \par $\gamma (a(bc))=f(a,bc)+\gamma
(a)\gamma (bc)=f(a,bc) +\gamma (a)(f(b,c)+\gamma (b)\gamma (c))$ and
\par $\gamma ((ab)c)=f(ab,c)+\gamma
(ab)\gamma (c)=f(ab,c) +(f(a,b)+\gamma (a)\gamma (b))\gamma (c)$,
\par consequently, \par $0={\sf t}_3(a,b,c)\gamma (a(bc))- \gamma ((ab)c) =
{\sf t}_3(a,b,c)f(a,bc) +{\sf t}_3(a,b,c)\gamma (a)(f(b,c)+\gamma
(b)\gamma (c))- f(ab,c) - (f(a,b)+\gamma (a)\gamma (b))\gamma (c)$.
\\ Taking into account that $\gamma (a)m=a\cdot m$ and $m\gamma (a)=m\cdot a$
for each $m\in M$ and $a\in A$ we deduce using formula 13$(1)$ that
\par $0= {\sf t}_3(a,b,c)\cdot a\cdot f(b,c)- f(ab,c)+{\sf
t}_3(a,b,c)\cdot f(a,bc)-f(a,b)\cdot c= (\delta ^2f)(a,b,c)$. \\
Thus $f\in B^2(A,M)$ and hence $f=(\delta ^1h)$ with $h\in
C^1(A,M):=Hom_{\cal T}(A,M)$. \par It remains to prove that the set $S$ of
all elements $\gamma (a)+h(a)$ forms in $P$ a subalgebra isomorphic
with $A$. From the construction of $S$ it follows that $S$ is a
a two-sided ${\cal T}$-module. We verify that it is closed relative to the
multiplication for all $a$ and $b$ in $A$:
\par $(\gamma (a)+h(a))(\gamma (b)+h(b))= \gamma (a)\gamma (b)
+\gamma (a)h(b)+h(a)\gamma (b)=\gamma (ab)-f(a,b)+ah(b)+h(a)b=\gamma
(ab)+h(ab)+ah(b)-h(ab)+h(a)b-f(a,b)= \gamma (ab)+h(ab)+(\delta
^1h)(a,b)-f(a,b)=\gamma (ab)+h(ab)$.
\par If there are given $A$, $M$ and $f$, then an enlargement
$P$ can be constructed as the direct sum $P=M\oplus A$ of two-sided ${\cal T}$-modules
and with the multiplication rule:
\par $(m_1+b_1)(m_2+b_2)=m_1b_2+m_2b_1+f(b_1,b_2)+b_1b_2$
\\ for every $m_1$ and $m_2$ in $M$ and $b_1$ and $b_2$ in $A$.
It rest to verify that this multiplication rule is homogeneous over
$\cal T$ and right and left distributive. Evidently,
\par $(m_1+b_1)(s(m_2+b_2))=(s(m_1+b_1))(m_2+b_2)=
s((m_1+b_1)(m_2+b_2)) $\par $=sm_1b_2+sm_2b_1+sf(b_1,b_2)+sb_1b_2$
and
\par $(sp)(m_1+b_1)=s(p(m_1+b_1))$
\\ for all $s, p\in \cal T$ and $m_1$ and $m_2$ in $M$ and $b_1$ and $b_2$ in
$A$, since ${\cal T}\subset {\cal C}(A)$ and $f(s,p)=0$. Moreover, \par
$(m_1+b_1)((m_2+b_2)+(m_3+b_3))=(m_1+b_1)((m_2+m_3)+(b_2+b_3))$\par
$=m_1(b_2+b_3)+(m_2+m_3)b_1+f(b_1,b_2+b_3)+b_1(b_2+b_3)$
\par $=m_1b_2+m_1b_3+m_2b_1+m_3b_1+f(b_1,b_2)+f(b_1,b_3)+b_1b_2+b_1b_3$
\par $=(m_1+b_1)(m_2+b_2)+(m_1+b_1)(m_3+b_3)$ and
analogously \par $((m_1+b_1)+(m_2+b_2))(m_3+b_3)=
(m_1+b_1)(m_3+b_3)+(m_2+b_2)(m_3+b_3)$ for all $m_1, m_2$ and $m_3$
in $M$ and $b_1, b_2$ and $b_3$ in $A$.

\par {\bf 16. Definition.} Let $M$ and $P$ and $N$ be two-sided $A$-modules, where $A$ is
a nonassociative metagroup algebra over a commutative associative
unital ring ${\cal T}$. A homomorphism (isomorphism)
$f: M\to P$ is called a right (operator) homomorphism (isomorphism)
if it is such for $M$ and $N$ as right $A$-modules, that is
$f(x+y)=f(x)+f(y)$ and $f(xa)=f(x)a$ for each $x$ and $y$ in $M$ and
$a\in A$. An enlargement $(P,\eta )$ of $M$ by $N$ is called right
inessential, if a right isomorphism $\gamma : N\to P$ exists
satisfying the restriction $\eta \gamma |_N=1|_N$.

\par {\bf 17. Theorem.} {\it Suppose that $M$ is a two-sided $A$-module,
where $A$ is a nonassociative metagroup algebra over a commutative associative
unital ring ${\cal T}$. Then for each
$n\ge 0$ there exists a two-sided $A$-module $P_n$ such that $H^{n+1}(A,M)$
is isomorphic with the additive group of equivalence classes of
right inessential enlargements of $M$ by $P_n$.}
\par {\bf Proof.} Consider two right inessential enlargements
$(E_1,\eta _1)$ and $(E_2,\eta _2)$ of $M$ by $N$, where $\xi _1$
and $\xi _2$ are embeddings of $M$ into $E_1$ and $E_2$
correspondingly. Take a submodule $Q$ of $E_1\oplus E_2$ consisting
of all elements $(x_1,x_2)$ satisfying the condition: $\eta
_1(x_1)=\eta _2(x_2)$. Then a quotient module $Q/T$ exists, where
$T= \{ (\xi m,-\xi m): m\in M \} $. Therefore $(\xi _1M\oplus \xi
_2M)/T$ is isomorphic with $M$ and homomorphisms $\eta _1$ and $\eta
_2$ induce a homomorphism $\eta $ of $Q/T$ onto $N$. Hence the
submodule $ker (\eta )$ is isomorphic with $M$. Then an addition of
enlargements is prescribed by the formula: $(E_1,\eta _1)+(E_2,\eta
_2):=(Q/T, \eta )$. Evidently sums of equivalent enlargements are
equivalent.
\par For an enlargement $(E, \eta )$ of $M$ by $N$ one takes
the direct sum of modules $E\oplus M$ and puts $T_b$ to be its
submodule consisting of all elements $(\xi m, - b \xi m)$ with $m\in
M$, where $\xi $ is an embedding of $M$ into $E$, $ ~ b\in \cal T$.
Therefore, a homomorphism $\eta $ induces a homomorphism
$\mbox{}_b\eta $ of $(E\oplus M)/T_b$ onto $N$, since the mapping
$(\xi m,m)\mapsto b\xi m+m$ is a homomorphism of $(\xi M)\oplus M$
onto $M$, also the ring $\cal T$ is commutative and associative. This induces an
enlargement of $M$ by $N$ denoted by $(\mbox{}_bE, \mbox{}_b\eta )$
and hence an operation of scalar multiplication of an enlargement
$(E,\eta )$ on $b\in \cal T$. From this construction it follows that
equivalent enlargements have equivalent scalar multiplies on $b\in
\cal T$.
\par Let $P_n$ be a $\cal T$-linear span of all elements $(x_1,...,x_{n+1})$ with
$x_1$,...,$x_{n+1}$ in $G$ such that
$((bx_1),x_2,...,x_{n+1})=(x_1,...,(bx_{n+1}))$ for each $b\in \cal T$. Next we
put
\\ $(1)$ $(x_1,...,x_{n+1})\cdot y := t_{n+2}(x_1,...,x_{n+1},y;l(n+2),u_{n+2}(n+2))  \cdot
(x_1,...,x_{n},(x_{n+1}y))$
and
\par $(2)$ $y\cdot (x_1,...,x_{n+1}) =\sum_{j=1}^{n+1}
(-1)^{j+1}\cdot $\par
$t_{n+2}(y,x_1,...,x_{n+1};u_1(n+2),u_{j+1}(n+2)) \cdot <y,x_1,
x_2,..., x_{n+1}>_{j,n+2}$\\  (see also notation 9$(2-4)$ and 12$(10-13)$) for every
$y, x_1,...,x_{n+1}$ in $G$. That is
$P_n$ is the two-sided $A$-module, where $A$ has the unit element.
\par We denote by $R_n=R(P_n,M)$ the
family of all right homomorphisms of $P_n$ into $M$. For each $p\in
R_n$ let an arbitrary element $\dot{p}\in C^{n}(A,M)$ in the additive group of all $n$ cochains
(that is, $n$ times $\cal T$-linear mappings of $A$ into $M$)
on $A$ with values in $M$ be prescribed by the formula $\dot{p} (a_1,...,a_n) = p(a_1,...,a_n,1)$ for all
$a_1,...,a_n$ in $A$, consequently, $(\dot{p} (a_1,...,a_n))\cdot y=
p(a_1,...,a_n,y)$ for each $y\in A$, since
$t_{n+3}(x_1,...,x_{n+1},1,g;l(n+3),u_{n+3}(n+3))=1$ for all
$x_1,...,x_{n+1}$ and $g$ in $G$. This makes the mapping $p\mapsto
\dot{p}$ an $\cal T$-linear isomorphism of $R_n$ onto $C^{n}(A,M)$.
\par Supply $C^n(A,M)$ with a two-sided $A$-module structure:
\par $(3)$ $(x_0\cdot f)(x_1,...,x_n)= x_0\cdot (f(x_1,...,x_n))$ and
\par $(4)$ $(f\cdot x_0)(x_1,..,x_n) = \sum_{k=0}^{n-1} (-1)^k t_{n+1}
(x_0,x_1,...,x_n; u_1(n+1),u_{k+2}(n+1))\cdot
f(x_0,...,x_kx_{k+1},...,x_n)$ $+(-1)^n (f(x_0,...,x_{n-1}))\cdot
x_n$\\ for each $f\in C^n(A,M)$ and all $x_0, x_1,...,x_n$ in $G$
extending $f$ by $\cal T$-linearity on $A$ from $G$, where $u_j(n+1)$ are given by formulas
12$(10-13)$. Thus the mapping
$p\mapsto \dot{p}$ is an operator isomorphism, consequently,
$H^l(A,R_n)$ is isomorphic with $H^l(A,C^{n}(A,M))$ for each
integers $n$ and $l$ such that $n\ge 0$ and $l\ge 0$. On the other
hand, $H^l(A,C^n(A,M))$ is isomorphic with $H^{l+n}(A,M)$ for each
$l\ge 1$, hence $H^l(A,R_n)$ is isomorphic with $H^{l+n}(A,M)$. In
virtue of Theorem 14 applied with $l=1$ we infer that $H^{n+1}(A,M)$ is
isomorphic with the additive group of equivalence classes of right
inessential enlargements of $M$ by $P_n$.

\par {\bf 18. Theorem.} {\it Let $M$ be a two-sided $A$-module,
where $A$ is a nonassociative metagroup algebra over a commutative
associative unital ring $\cal T$. Then to each
$n+1$-cocycle $f\in Z^{n+1}(A,M)$ an enlargement of $M$ by a two-sided
$A$-module $P_n$ corresponds such that $f$ becomes a
coboundary in it.}
\par {\bf Proof.} A $n+1$-cocycle $f\in Z^{n+1}(A,M)$ induces an enlargement $(E,\eta )$
of $M$ by $P_n$ due to Theorem 17. An element $h$ in $Z^1(A,R_n)$
corresponding to $f$ is characterized by the equality:
\par $(h(x_1))(x_2,...,x_{n+1},1) = t_{n+1}(x_1,...,x_{n+1};u_1(n+1),l(n+1))\cdot
f(x_1,...,x_{n+1})$ \\ for all $x_1,...,x_{n+1}$ in $G$. This enlargement
$(E,\eta )$ as the two-sided $A$-module is $P_n\oplus M$ such that
\par $x_1\cdot ((x_2,...,x_{n+1},1),0)= (x_1\cdot
(x_2,...,x_{n+1}),f(x_1,...,x_{n+1}))$.
\par Let $\gamma (a_1,...,a_n)=(a_1,...,a_n,0)$ for all
$a_1,...,a_n$ in $A$. Therefore we deduce that
\par $f(x_1,...,x_{n+1}) = t_{n+1}(x_1,...,x_{n+1};l(n+1),u_1(n+1))\cdot
\{ x_1\cdot \gamma (x_2,...,x_{n+1},1)$\par $ - \gamma (x_1\cdot
(x_2,...,x_{n+1},1)) \} $. \\ There exists a $n$-cochain $v\in
C^{n}(A,E)$ defined by $v(a_1,...,a_{n})=
((a_1,...,a_{n},1),0)$ for all $a_1,...,a_{n}$ in $A$. Thus
$f=\delta v$.

\par {\bf 19. Theorem.} {\it Let $A$ be a nonassociative metagroup
algebra over a commutative associative unital ring $\cal T$.
Then an algebra $B$ over $\cal T$ exists such that $B$
contains $A$ and each $\cal T$-homogeneous derivation $d: A\to A$ is the
restriction of an inner derivation of $B$.}
\par {\bf Proof.} Naturally an algebra $A$ has the structure of a two-sided
$A$-module. In view of Proposition 12 each derivation of the two-sided
algebra $A$ can be considered as an element of $Z^1(A,A)$.
\par Applying Theorem 18 by induction one obtains a two-sided
$A$-module $Q$ containing $M$ for which an arbitrary element
of $Z^{n+1}(A,M)$ is represented as the coboundary of an element of
$C^n(A,Q)$. At the same time $M$ and $Q$ satisfy conditions 7$(1-3)$.
This implies that the natural injection of $H^{n+1}(A,M)$ into
$H^{n+1}(A,Q)$ maps $H^{n+1}(A,M)$ into zero.
\par Therefore a two-sided $A$-module $E$ exists which as a two-sided ${\cal T}$-module
is a direct sum $A\oplus P$ and $P$ is such
that for each $f\in Z^1(A,A)$ there exists an element $p\in P$
generally depending on $f$ with the property $f(a)=a\cdot p-p\cdot
a$. To the algebra $A$ the metagroup $G$ corresponds. Enlarging $P$
if necessary we can consider that to $P$ a metagroup $G$ also
corresponds in such a manner that properties 7$(1-3)$ are fulfilled. Now
one can take $A\oplus P$ as the underlying two-sided $\cal T$-module of $B$ and
supply it with the multiplication $(a_1,p_1)(a_2,p_2) := (a_1a_2,
a_1\cdot p_2+p_1\cdot a_2)$ as the semidirect product for each
$a_1$, $a_2$ in $A$ and $p_1$, $p_2$ in $P$. An embedding $\xi $ of
$A$ into $B$ is $\xi (a)=(a,0)$ for each $a$ in $A$. This implies
that $f(a)=(a,0)(0,p)-(0,p)(a,0)= a(0,p)-(0,p)a$.

\par {\bf 20. Theorem.} {\it Suppose that $A$ is a nonassociative metagroup
algebra of finite order over a commutative associative unital ring ${\cal T}$
and $M$ is a finitely generated two-sided $A$-module.
Then $M$ is semisimple if and only if its cohomology group
is null $H^n(A,M)=0$ for each natural number $n\ge 1$.}
\par {\bf Proof.} Certainly, if $E$ is an $A$-module and $N$ its
$A$-submodule, then a natural quotient morphism $\pi : E\to E/N$
exists. Therefore, an enlargement $(E,\eta )$ of a two-sided $A$-module
$M$ by a two-sided $A$-module $N$ is inessential if and only if there
is a submodule $T$ in $E$ complemented to $\xi (M)$ such that $T$ is
isomorphic with $E/\xi (M)$, where $\xi $ is an embedding of $M$
into $E$. When $M$ is semisimple, it is either simple or a finite
product of simple modules, since $M$ is finitely generated. For a
finitely generated module $E$ and its submodule $N$ the quotient
module $E/N$ is not isomorphic with $E$, since the algebra $A$ is of
finite order over the commutative associative unital ring $\cal T$.
\par In virtue of Theorems 17 and 18 if for an algebra $A$ its corresponding
finitely generated two-sided $A$-modules are semisimple, then its
cohomology groups of dimension $n\ge 1$ are zero. \par Vise versa
suppose that $H^n(A,M)=0$ for each natural number $n\ge 1$. Consider
a finitely generated two-sided $A$-module $E$ and its two-sided
$A$-submodule $N$. At first we take into account the right
$A$-module structure $E_r$ of $E$ with the same right
transformations, but with zero left transformations. Then the left
inessential $(E_r,\eta _r)$ enlargement of $M_r$ by $N_r=E_r/\xi
_r(M_r)$ exists, where $\eta _r: E_r\to E_r/\xi _r(M_r)$ is the
quotient mapping and $\xi _r$ is an embedding of $M_r$ into $E_r$.
From Theorem 17 it follows that the enlargement $(E_r,\eta _r)$ is
right inessential. Analogously considering left $A$-module
structures $E_l$ and $M_l$ we infer that $(E_l,\eta _l)$ is also
left inessential.

\par {\bf 21. Note.} Let $A$ be a nonassociative metagroup algebra over
a commutative associative unital ring $\cal T$ of characteristic
$char ({\cal T})$ other than two and three. There exists its opposite algebra $A^{op}$. The latter as an
$F$-linear space is the same, but with the multiplication $x\circ y
=yx$ for each $x, y \in A^{op}$. To each element $h\in A$ or $y\in
A^{op}$ pose a left multiplication operator $L_h$ by the formula
$xL_h=hx$ or a right multiplication operator $xR_y=xy$ for each
$x\in A$ respectively. Having the anti-isomorphism operator $S: A\to
A^{op}$, $A\ni x\mapsto xS\in A^{op}$, $S(xy)=S(y)S(x)$, one gets
\par $(1)$ $R_hS=SL_{hS}$ and $L_hS=SR_{hS}$ \\ for each $x$ and $h$
in $A$. Then taking into account $(1)$ analogously to formula
17$(4)$ we put
\par $x_0\cdot (L_{x_1},R_{x_2}) =
{\sf t}_3^{-1} (x_0,x_1,x_2)\cdot (x_0L_{x_1})L_{x_2S} -
x_0(L_{x_1}R_{x_2})$\par $+ {\sf t}_3^{-1} (x_0,x_1,x_2)\cdot
(x_0R_{x_1S})R_{x_2} $ \\
and taking into account multipliers ${\sf t_3}$ this gives
\par $(2)$ $x_0\cdot (L_{x_1},R_{x_2}) = x_0(L_{x_1}L_{x_2S} - L_{x_1}R_{x_2} + R_{x_1S}R_{x_2})$
\\ for all $x_0, x_1,x_2$ in $G$. Then symmetrically $S(x_0\cdot (L_{x_1},R_{x_2}))$
provides the formula for $(L_{y_1},R_{y_2})\cdot y_0$ for each $y_0,
y_1$ and $y_2$ in $G$.  Extending these rules by $\cal T$-linearity on
$A$ and $A\otimes _{\cal T}A^{op}$ from $G$ one supplies the tensor product
$M=A^e$ over $\cal T$, where $A^e=A\otimes_{\cal T} A^{op}$ is the
enveloping algebra, with the two-sided $A$-module structure.

\par {\bf 22. Corollary.} {\it Let $A$ be a semisimple nonassociative metagroup
algebra of finite order over a commutative associative unital ring
$\cal T$ of characteristic $char ({\cal T})$ other than two and
three and let $M$ be a two-sided $A$-module described in \S 21. Then
$H^n(A,M)=0$ for each natural number $n\ge 1$.}
\par {\bf Proof.} Since $A$ is semisimple, then the module $M$ from
\S 21 is semisimple, consequently, from Theorem 20 the statement of
this corollary follows.

\par {\bf 23. Theorem.} {\it Let $A$ be a nonassociative metagroup algebra
over a commutative associative unital ring $\cal T$
with a metagroup $G$ such that $G\cap {\cal T} = F^{\times }$,
and let $H$ be a proper normal submetagroup isomorphic with $G$ such that the
quotient metagroup $G/H$ is infinite and $F^{\times }\subset
H$. Then outer $\cal T$-linear automorphisms of $A$
exist and the cardinality of their family is not less, than
$2^{\xi }$, where $\xi =card (G/F^{\times })$.}
\par {\bf Proof.} By the conditions of this theorem
$G/H$ is infinite, that is the cardinality
$card (G/H) \ge \aleph _0$. Moreover we infer that $G/H=(G/F^{\times
})/(H/F^{\times })$, since $F^{\times }\subset H\subset G$
and $H$ is isomorphic with $G$. At the same time the quotients
$J:=G/F^{\times }$ and $J_{\infty }:=H/F^{\times }$ are
groups. Thus $G/H$ is the group. \par If an element $b$ of
a group $P$ is infinite divisible, that is for each $n\in \bf N$
there exists $a\in P$ having the property $a^n=b$, we put $ord (b)
=0$. Then $ord (b) =m>0$ is a least natural number $m$ such that
$b^m=e$ is a neutral element. Otherwise $ord (b)=\infty $, when $b$
is neither infinite divisible nor is of finite order. From the
theorem about transfinite induction \cite{kunenb} it follows that
there are sets $\Upsilon $ and $\Upsilon _{\infty }$ such that
\par $(1)$ $\Upsilon \subset J$ and $\Upsilon _{\infty }\subset J_{\infty }$
satisfying the properties:
\par $(2)$ $card (\Upsilon )\ge \aleph _0$ and $card (\Upsilon _{\infty })\ge \aleph
_0$, $card (\Upsilon )= card (\Upsilon _{\infty })$, also \par $card
(J\setminus gr _J(\Upsilon ))\le card (J)$;
\par $(3)$ $\Upsilon = \bigcup_{b\in \Psi } b\Upsilon_{\infty }$,
where $\Psi $ is a set of the cardinality $card (\Psi )\ge \aleph _0$, \par $b\Upsilon _{\infty }\cap
c\Upsilon _{\infty }=\emptyset $ for each $b\ne c$ in $\Psi $;
\par $(4)$ $\forall b\in \Upsilon $ $b\notin gr _J(\Upsilon \setminus
\{ b \} )$ and $\forall b\in \Upsilon _{\infty }$ $b\notin
gr_{J_{\infty }} (\Upsilon _{\infty }\setminus \{ b \} )$, \\ since
$card (B\times B)=card (B)$ for $card (B)\ge \aleph _0$, where $gr_J
(K)$ denotes an intersection of all subgroups in $J$ containing a
subset $K$, where $K\subset J$.
The cardinality of the family of inner
automorphisms of $G$ is not greater than $\xi $, since each inner automorphism has the form
$G\ni g\mapsto f^{-1}((h^{-1}(gh))f)\in G$, where $f$ and $h$ are fixed in $G$.
\par From conditions $(1)$ and $(2)$ we deduce that there is $n\in
\{ 0, 1, 2,..., \infty \} $ such that \par $(5)$ $card (\Upsilon
(n))=card (\Upsilon _{\infty }(n))=card (\Upsilon )$, where
$\Upsilon (n):= \{ b\in \Upsilon: ord (b) =n \} $. \par From Conditions
$(2-5)$ it follows that there are automorphisms of
$gr_J(\Upsilon (n))$ and of $gr_{J_{\infty }} (\Upsilon _{\infty
}(n))$ which are not inner and induced by bijective
surjective mappings of the sets $\Upsilon (n)$ and $\Upsilon _{\infty
}(n)$, since $card (gr_J(\Upsilon (n))) = card (gr_J(\Upsilon
))=card (J) =:\xi \ge \aleph _0$ and $card (Aut (gr_J(\Upsilon (n)))=2^{\xi }>\xi
$. In virtue of the theorem about extensions of
automorphisms \cite{neumann} (or see \cite{plotk,puusfocus}) these
automorphisms have extensions from the aforementioned subgroups
on groups isomorphic with $J$ and $J_{\infty }$ respectively.
\par This implies that nontrivial outer automorphisms $\theta '$ of $J$
and ${\theta _{\infty }}'$ of $J_{\infty }$ exist, since $card
(J\setminus gr _J(\Upsilon ))\le card (J)$. They induce outer
automorphisms $\theta $ of $G$ and $\theta _{\infty }$ of $H$
correspondingly, since $F^{\times }$ is the commutative
multiplicative normal subgroup in them. From their construction
it follows that their restrictions on $F^{\times }$ is the identity mapping
$\theta |_{F^{\times }}=id$ and $\theta _{\infty }|_{F^{\times }}=id$.
On the other hand, $G\cap {\cal T} = F^{\times }$, consequently,
there exist $\cal T$-linear extensions of automorphisms $\theta $ and $\theta _{\infty }$
from metagroups $G$ and $H$ on metagroup algebras ${\cal T}[G]$ and ${\cal T}[H]$
respectively.
\par The cardinality of the family of all outer $\cal T$-linear automorphisms of $A$
is not less, than $2^{\xi }$, since $2^{\xi }>\xi \ge \aleph _0$.

\par {\bf 24. Corollary.} {\it Let conditions of Theorem 23 be satisfied.
\par $(1)$. Moreover, suppose that $G$ is a topological Hausdorff metagroup and for $J$ in the quotient topology
$\Upsilon (n)$ fulfilling conditions 23$(1-5)$ is everywhere dense
in $J$. Then the metagroup $G$ has discontinuous outer
automorphisms.
\par $(2)$. If $A={\cal T}[G]$ is a topological Hausdorff metagroup algebra
such that conditions in $(1)$ are fulfilled for $G$ in the hereditary from $A$ topology,
then $A$ has discontinuous outer $\cal T$-linear automorphisms.}
\par {\bf Proof.} $(1)$. Take any two nonintersecting open subsets $U_1=U_1^{-1}$
and $U_2$ in $J$ with $U_2=gU_1$ for some $g\in J$ and subsets
$V_{i,j}\subset U_i$ such that $V_{i,1}\cap V_{i,2}=\emptyset $ and
$V_{i,1}\cup V_{i,2}\subseteq \Upsilon (n)\cap U_i$, $card (V_{1,1})
= card (V_{i,j})$ and $V_{i,j}$ is everywhere dense in $U_i$ for
each $i$ and $j$ in $ \{ 1, 2 \} $. Using conditions 28$(1-5)$
choose an automorphism $\theta '$ of $J$ having the property $\theta
'(V_{i,1}) = V_{i,1}$ for $i=1$ and $i=2$, $\theta '(V_{1,2}) =
V_{2,2}$ and $\theta '(V_{2,2})=V_{1,2}$. Thus $\theta '$ is
discontinuous. It induces a discontinuous automorphism $\theta $ of
$G$.
\par $(2)$ follows from $(1)$.

\end{document}